\documentclass[12pt]{amsart}
\title{On D.K. Biss' papers\\
``The homotopy type of the matroid Grassmannian" \\
 Annals of Mathematics 158 (2003) 929-952 \\ and \\
``Oriented matroids, complex manifolds, and a
combinatorial model for BU"\\ Advances in Mathematics 179 (2003) 250-290 \\
 }
\author{N. Mn\"ev}
\begin{document}
\maketitle

My very unfortunate duty is to point out a serious flaw
in  papers \cite{B,B1} by Daniel Biss devoted to homotopy type of matroid Grassmannians.
Four years passed after the publication, but no  errata became available.
Meanwhile the problem was already
acknowledged and  discussed in private correspondence
between experts in April 2006.

The mistake is the same for both papers,
it is very simple, looks almost like a typo,
but it is located in the key propositions
for the induction towers --  Proposition 4.5 \cite{B} on page 948
and Proposition 7.3 \cite{B1} on page 285.
Here the aim of Biss' reasons is to show that his natural
combinatorial models for Schubert cells
are correctly attached one to another.
In the surgery of such an attachment
some important complex called  $|| S^+ ||$ appears.
Biss needs to prove that  $|| S^+ ||$ is contractible.
Biss  covers $|| S^+ ||$ by
contractible {\em open} set
$O=(|| S^+ || \setminus || A ||)$ and
a {\em closed} set $C=||B||$ homotopy  equivalent to $C \cap O$
and concludes
that it follows that $||S^+ ||$ is contractible.
This is not correct. For example one can
easily cover
a circle by  open and closed intervals with the same property and thus prove that
the circle is also contractible.

So the statement in the final lines of the proof of \cite[Proposition 4.5, 11-13 lines from the top of the page 948]{B}
:
\begin{quotation}\em
``...and thus
the inclusion $||S^+||\setminus ||A||\hookrightarrow ||S^+||$
induces a homotopy equivalence, and $ ||S^+|| $ is therefore contractible''
\end{quotation}
is wrong.

The same argument is used by Biss in the proof of
\cite[Proposition 7.3, page 285 lines 11-13 from the top ]{B1}.
There the role of $||S^+||$ is played by $||S_+^\mathfrak{R}||$.

Unfortunately this simple mistake destroys the main theorems
of both papers. If one tries to continue the cut and paste induction over $||S^+||$ by
cutting it into smaller natural pieces then one will quickly face "ball-like"
posets formed by {\em weak maps} of oriented matroids which have a well known bad
habit  to be homotopy nontrivial (see for example \cite{Mnev:Richter}).
Surprisingly the weak maps of matroids do not appear up to this hidden point in
Biss' scheme at all.  This looks strange by the reasons of dimension
of Grassmannians.

Personally, I  don't currently see a way to save Biss' theorems in those strongest
forms (that finite matroid posets can model finite-dimensional
Grassmannians).  On the other hand in stable infinite dimensional case when Grassmannian became
classifying space for stable  vector bundles I think that matroid models really may
work very well and one can use much more geometric topology  to prove
it.

\def\cprime{$'$} \def\cprime{$'$}

\end{document}